\documentclass{article}

\usepackage{amsfonts}
\usepackage{amsmath}
\usepackage{amsthm}

\newtheorem{lemma}{Lemma}[section]
\newtheorem{theorem}[lemma]{Theorem}
\newtheorem{definition}[lemma]{Definition}
\newtheorem{corollary}[lemma]{Corollary}
\newtheorem{observation}[lemma]{Observation}

\begin{document}
\title{Conditions for Weighted Cover Pebbling of
Graphs}
\author{Annalies Vuong\footnote{Department of
Mathematics, U.C. Santa Barbara,
Santa Barbara, CA 93106. email:
azv@umail.ucsb.edu. Supported by
National Science Foundation grant \#DMS-0139286}
\and M. Ian
Wyckoff\footnote{Department of Mathematics and
Statistics, Arizona
State University, Tempe, AZ 85287-1804. email:
iwyckoff@yahoo.com.
Supported by National Science Foundation grant
\#DMS-0139286}}
\date{\today}
\maketitle

\begin{abstract}
In a graph $G$ with a distribution of pebbles on
its vertices, a
pebbling move is the removal of two pebbles from
one vertex and
the addition of one pebble to an adjacent vertex.
 A weight
function on $G$ is a non-negative integer-valued
function on the
vertices of $G$.  A distribution of pebbles on
$G$ covers a weight
function if there exists a sequence of pebbling
moves that gives a
new distribution in which every vertex has at
least as many
pebbles as its weight. In this paper we give some
necessary and
some sufficient conditions for a distribution of
pebbles to cover
a given weight function on a connected graph $G$.
As a corollary,
we give a simple formulation for the `weighted
cover pebbling
number' of a weight function $W$ and a connected
graph $G$,
defined by Crull et al. to be the smallest number
$m$ such that
any distribution on $G$ of $m$ pebbles is a cover
for $W$. Also,
we prove a cover pebbling variant of Graham's
Conjecture for
pebbling.

\end{abstract}

\section{Introduction}
Suppose $k$ pebbles are placed on the vertices of
a graph $G$.
Define a pebbling move on this distribution as
the removal of two
pebbles from one vertex together with the
addition of one pebble
to an adjacent vertex. In the game of pebbling,
one attempts to
place a pebble on a specified vertex in a graph
by a sequence of
pebbling moves on some starting distribution of
pebbles. This game
was first suggested by Lagarias and Saks, in an
attempt to provide
a short alternate proof of a result in additive
number theory.
Chung~\cite{ch} later succeeded in creating such
a proof and in so
doing defined the pebbling number $\pi(G)$ of a
connected graph
$G$, the minimum number $m$ such that one can win
the pebbling
game with any target vertex and any starting
distribution of $m$
pebbles. See \cite{survey} for a more complete
treatment of
pebbling.

Crull et al.~\cite{hu} introduced weighted cover
pebbling, a
variation on the game of pebbling.  In weighted
cover pebbling, we
begin with a connected graph $G$ and a
non-negative integer valued
`weight function' on the vertices of $G$. 
Starting from an
initial distribution of pebbles, a player
performs a sequence of
pebbling moves in an attempt to create a final
distribution on $G$
such that every vertex has at least as many
pebbles as its weight.
Crull et al. also defined the weighted cover
pebbling number of a
connected graph $G$, $\Phi_W(G)$, as the minimum
number $m$ such
that one can win the cover pebbling game
(relative to the weight
function $W$) with any starting distribution of
$m$ pebbles.

In this paper we provide some necessary and some
sufficient
conditions for a given distribution of pebbles to
allow a solution
for the weighted cover pebbling game, relative to
a strictly
positive weight function.  As a direct corollary
of these results,
we give a simple formulation of the value
$\Phi_W(G)$ for any
strictly positive weight function $W$, which
vastly simplifies
earlier proofs for particular weighted cover
pebbling numbers of
hypercubes~\cite{cp-hypercubes}, complete
multipartite
graphs~\cite{cp-results}, and trees~\cite{hu}.
This corollary is
equivalent to an affirmative answer to question
10 posed in
\cite{hu}.

\section{Preliminaries}

We now give a list of notations and definitions
that characterize
the game of cover pebbling.  Our characterization
attempts to
generalize and formalize many of the standard
terms previously
developed in the literature, both for pebbling
and cover pebbling.

\begin{definition}
Let $G$ be a graph, and let $\mathfrak{I}$ denote
the set of all
integer-valued functions on $V(G)$.  Let
$\mathfrak{D}$ denote all
functions in $\mathfrak{I}$ that are nonnegative.
 We call a
function $D$ a \emph{distribution} on $G$ if and
only if $D \in
\mathfrak{D}$.
\end{definition}

We may also call a distribution a \emph{weight
function}.
Following this, we will refer to a current or
previous
configuration of pebbles on a graph as a
distribution, and reserve
use of the term `weight function' for a desired
distribution to be
reached after completion of pebbling moves.

\begin{definition}
Let $G$ be a graph, and let $p, q$ be two
adjacent vertices of
$G$. For all $D \in \mathfrak{D}$ , define the
function $P_{p,q}:
\mathfrak{D} \mapsto \mathfrak{I}$ by

\[
[P_{p,q}(D)](s) = \left\{ \begin{array}{ll}
D(s)-2 & s = p \\
D(s) + 1 & s = q \\
D(s) & \rm{otherwise} \end{array} \right.
\]
If $P(D) \in \mathfrak{D}$, then $P$ is called a
\emph{pebbling
move} on $D$. We also consider the identity
function to be a
pebbling move on $D$ for all $D \in
\mathfrak{D}$.
\end{definition}

\begin{definition}\label{derived}
Let $G$ be a graph, and let $D$, $D'$ be
distributions on $G$.  If
there exists a sequence of functions $P_1,
\ldots, P_n$ such that
$P_i$ is a pebbling move on $P_{i-1} \circ \ldots
\circ P_1(D)
\;\forall i \in [n],$ and
\[P_{n} \circ \ldots \circ P_1(D)= D',\] then we
say that $D'$ is derivable from $D$\!.
\end{definition}

\begin{definition}
Let $G$ be a graph, and let $D$, $D'$ be
distributions on $G$.  If
$D(q) \leq D'(q) \; \forall q \in V(G)$, then we
say that $D$ is
contained in $D'$, and if $q$ is a vertex such
that $D(q)
> D'(q)$, then we say that $q$ is a distribution
node of $D$
relative to $D'$.  Also, if there exists a
distribution derivable
from $D$ that contains $D'$, then we say that $D$
is a cover of
$D'$, or $D$ covers $D'$.
\end{definition}

\begin{definition}
Let $H$ be an induced subgraph of $G$, $G$ an
induced subgraph of
$K$, and $D$ a distribution on $G$. By $D_H$ we
denote the
distribution on $H$ defined by
\[ D_H(q) = D(q) \quad \forall q \in V(H). \]
By $D_K$ we denote the distribution on $K$
defined by
\[D_K(q) = \left\{ \begin{array}{ll} D(q) & q \in
V(G) \\
0 & q \in V(K) \setminus V(G) \end{array} \right.
.\] Also, if $S
\subseteq V$, by $D_S$ we denote the distribution
$D_{G[S]}$ on
$G[S]$ and by $D\rvert_S$, the distribution on
$G$ such that
\[ D\rvert_S(q) = \left\{ \begin{array}{ll} D(q)
& q \in S \\
0 & q \in S^c \end{array} \right. .\]
\end{definition}

\begin{definition}
Let $S$ be a nonempty subset of $V(G)$, and let
$D$ be a
distribution on $G$.  Then we define the
\emph{standard value} of
$D$ with respect to $S$ as
$$ V_S(D) =  \sum_{q \in V(G)} D(q) \cdot 2^{d(q,
S)}$$
where $$d(q,S) = \min_{r \in S}d(q,r).$$ We also
define the
\emph{stacking number} of G relative to $W$ to be
$$ SN_W (G) = \max_{q \in V(G)} V_{\{q\}}(W) =
\max_{q \in V(G)}\sum_{u \in
V(G)}D(u)\cdot2^{d(u,q)}.$$
\end{definition}

\begin{observation} \label{value}
Let $S_1,S_2$ be two nonempty subsets of $V(G)$,
and let $D_1,D_2$
be two distributions on $G$.  Then the following
statements are
true:
\begin{enumerate}
\item If $S_1 \subset S_2$, then $V_{S_1}(D_1)
\geq V_{S_2}(D_1)$.
\item If $D_1$ is properly contained in $D_2$,
then $V_{S_1}(D_1)
< V_{S_1}(D_2)$. \item If there exists a legal
pebbling move $P$
on $D_1$ s.t. $P(D_1) = D_2$, then $V_{S_1}(D_1)
\geq
V_{S_1}(D_2)$. \item If $D_1$ is a cover of
$D_2$, then
$V_{S_1}(D_1) \geq V_{S_1}(D_2)$.
\end{enumerate}
\end{observation}

Note that, by property $4$, a distribution $D$
covers a weight
function $W$ only if $$V_{S}(D) \geq
V_{S}(W)\quad \forall\;S
\subseteq V(G),\,S \neq \emptyset.$$

\section{Principal Result}
In this section we present our primary result, a
strong sufficient
condition for a distribution $D$ to cover a given
positive weight
function $W$.

\begin{lemma}\label{nonode}
Let $G$ be a connected graph, $W$\,a weight
function on $G$, and
$D$ a distribution on $G$.  If $D$ has no
distribution nodes
relative to $W$, then $D$ is a cover of $W$ if
and only if $D =
W$.
\end{lemma}
\begin{proof}  By assumption, $D$ has no
distribution nodes relative to
$W$, thus $D(q) \leq W(q) \;\forall q \in V(G)$. 
By
Observation~\ref{value} we have that $D$ covers
$W$ only if
$V_{V(G)}(D) \geq V_{V(G)}(W)$.  Thus $D$ covers
$W$ only if $D =
W$.  The reverse implication is obvious.
\end{proof}

\begin{lemma} \label{Timothy}
Let $G$ be a connected graph, $W$\,a positive
weight function on
$G$, $v_0 \in V(G)$, and $D$ a distribution on
$G$ such that $D(s)
\leq W(s) \quad \forall s \neq v_0$.  If
$V_{\{v_0\}}(D) \leq
V_{\{v_0\}}(W)$, then there exists a distribution
$D^*$ on $G$
contained in $W$ \,and derivable from $D$ such
that
$V_{\{v_0\}}(D^*) = V_{\{v_0\}}(D)$.
\end{lemma}
\begin{proof}
By induction on $D(v_0)$.\\
If $D(v_0) \leq W(v_0),$ then we simply set $D^*
= D$.\\
Now assume that $D(v_0) > W(v_0)$, and that the
lemma holds for
any suitable distributions $D'$ with $D'(v_0) <
D(v_0)$. Select
$q$ in $V(G)$ such that $D(q) < W(q)$ and $d(v_0,
q)$ is
minimized. Such a point must exist, for otherwise
$W$ is properly
contained in $D$ and $V_{\{v_0\}}(D) >
V_{\{v_0\}}(W)$. Let $v_0,
\ldots, v_k, q$ be a shortest path from $v_0$ to
$q$. As $q$ is at
minimum distance to $v_0$, we have that $D(v_i) =
W(v_i) \;\forall
i \in [k]$. Define the distribution $D^I$ on $G$
by
\[D^I(s) = \left\{ \begin{array}{ll} D(s) - 2 & s
= v_0 \\
D(s)-1 & s \in \{v_1, \ldots, v_k\} \\
D(s) + 1 & s = q \\ D(s) & \rm{otherwise}
\end{array} \right. .\]
$D^I$ is derivable from $D$ by $P_{v_k, q} \circ
P_{v_{k-1}, v_k}
\circ \cdots \circ P_{v_0, v_1} (D) = D^I$. Also,
$$V_{\{v_0\}}(D^I) = V_{\{v_0\}}(D) - [2 + 2^1 +
\cdots + 2^{d(v_0, q) - 1} -
2^{d(v_0, q)}] = V_{\{v_0\}}(D)$$ and $D^I(s)
\leq W(s) \; \forall
s \neq v_0$, thus $D^I$ satisfies the induction
hypothesis for
$D(v_0) - 2$. Now we have a distribution $D^*$ on
$G$ contained in
$W$ and derivable from $D^I$ such that
$V_{\{v_0\}}(D^I) =
V_{\{v_0\}}(D^*)$. Clearly $D^*$ satisfies for
$D$.
\end{proof}

Note that the proof of Lemma~\ref{Timothy}
implies an algorithm
for achieving the desired distribution.

\begin{lemma}\label{onenode}
Let $G$ be a connected graph, $W$ a positive
weight function on
$G$, $v_0 \in V(G)$, and $D$ a distribution on
$G$ such that $D(s)
\leq W(s) \; \forall s \neq v_0$.  Then $D$
covers $W$ if and only
if $V_{\{v_0\}}(D) \geq V_{\{v_0\}}(W)$.
\end{lemma}
\begin{proof} By Observation~\ref{value}, we have
that $D$ covers $W$ only if
$V_{\{v_0\}}(D) \geq V_{\{v_0\}}(W)$.  Now,
assume that
$V_{\{v_0\}}(D) \geq V_{\{v_0\}}(W)$.  Consider
the distribution
$D^\#$ defined by
$$ D^\#(q) = \left\{ \begin{array}{ll}
D(v_0) - [V_{\{v_0\}}(D) - V_{\{v_0\}}(W)] & q =
v_0 \\
D(q) & \rm otherwise \end{array} \right. .$$
$D^\#$ satisfies the
conditions of Lemma~\ref{Timothy}, thus there
exists a
distribution $D^*$ on $G$ contained in $W$ and
derivable from
$D^\#$ such that $V_{\{v_0\}}(D^*) =
V_{\{v_0\}}(D^\#)$.  As
$V_{\{v_0\}}(D^*) = V_{\{v_0\}}(D^\#) =
V_{\{v_0\}}(D) -
[V_{\{v_0\}}(D) - V_{\{v_0\}}(W)] =
V_{\{v_0\}}(W)$ and $D^*$ is
contained in $W$, we have that $D^* = W$. $D^*$
is derivable from
$D^\#$, thus $D^\#$ is a cover of $W$, and as
$D^\#$ is contained
in $D$, this gives that $D$ is also a cover of
$W$.
\end{proof}

\begin{theorem}\label{bigproof}
Let $G$ be a connected graph, $W$\,a positive
weight function on
$G$, and $D$ a distribution on $G$ with a
nonempty set of
distribution nodes $\{d_1, \ldots, d_n\}$
relative to $W$. If
$V_{\{d_1,\dots,d_n\}}(D) \geq max_{i \in
[n]}V_{\{d_i\}}(W)$,
then $D$ covers $W$.
\end{theorem}
\begin{proof}  By induction on the number of
distribution nodes of $D$.\\
Suppose first that $D$ has only distribution node
$d_1$.  Then
$$ V_{\{d_1\}}(D) \geq V_{\{d_1\}}(W),$$
and by Lemma~\ref{onenode} $D$ covers $W$.\\
Suppose next that $D$ has distribution nodes
$d_1, \ldots, d_n$
with $n \geq 2$, and that the theorem holds for
any distribution
$\widetilde{D}$ having less than $n$ distribution
nodes relative
to $W$.  For all $i \in [n]$ define $E(d_i) = \{v
\in V(G) \,
\vert \, d(v, d_i) \leq d(v, d_j) \; \forall j
\in [n] \}$.  Let
$G_i = G[E(d_i)], D_i = D_{G_i}, W_i = D_{G_i} \,
\forall i \in
[n]$.  If $D_i$ covers $G_i \; \forall i \in
[n]$, then it is
clear that $D$ covers $G$.  Otherwise $\exists j
\in [n]$ such
that $D_j$ does not cover $G_j$.  Assume without
loss of
generality $j = n$.  By Lemma~\ref{Timothy},
there exists a
distribution $D^*$ on $G_n$ derivable from $D_n$
such that $D^*$
is contained in $W_n$ and
$$V_{\{d_n\}}(D^*) = V_{\{d_n\}}(D_n).$$
Thus, $$V_{\{d_n\}}(D^*_G) = V_{\{d_n\}}(D
\rvert_{E(d_n)}).$$
Consider the distribution
\[ D' = D \rvert_{E(d_n)^c} + D^*_G\] derivable
from $D$.  Clearly
$\{d_1, \ldots, d_{n-1}\}$ is the set of
distribution nodes for
$D'$. Also,
\begin{align*}
V_{\{d_1, \ldots, d_{n-1}\}}(D') & \geq V_{\{d_1,
\ldots,
d_n\}}(D') \\
& = V_{\{d_1, \ldots, d_{n}\}}(D
\rvert_{E(d_n)^c}) + V_{\{d_1,
\ldots,
d_n\}}(D^*_G)\\
& = V_{\{d_1, \ldots, d_{n}\}}(D
\rvert_{E(d_n)^c}) +
V_{\{d_n\}}(D^*_G)\\
& = V_{\{d_1, \ldots, d_{n}\}}(D
\rvert_{E(d_n)^c}) +
V_{\{d_n\}}(D \rvert_{E(d_n)})\\
& = V_{\{d_1, \ldots, d_{n}\}}(D
\rvert_{E(d_n)^c}) + V_{\{d_1,
\ldots, d_n\}}(D \rvert_{E(d_n)})\\
& = V_{\{d_1, \ldots, d_{n}\}}(D)\\
& \geq \max_{i \in [n]} V_{\{d_i\}}(W) \\
& \geq \max_{i \in [n-1]} V_{\{d_i\}}(W).
\end{align*}
Thus $D'$ satisfies the induction hypothesis for
$n-1$ and $D'$
covers $W$.  Since $D'$ is derivable from $D$,
$D$ also covers
$W$.
\end{proof}

\section{Corollaries}

\begin{corollary}{Stacking Theorem}\\
If $G$ is a connected graph and $W$\,is a
positive weight function
on $G$, then $\Phi_W(G) = SN_W(G)$.
\end{corollary}
\begin{proof}
Let $D$ be a distribution on $G$ such that $\vert
D \vert \geq
SN_W(G)$, where we define $\vert D \vert =
\sum_{q \in V(G)}
D(q)$. If $D$ has no distribution nodes relative
to $W$, then we
have that $\vert D \vert \geq SN_W(G) \geq \vert
W \vert$, thus by
Lemma~\ref{nonode} we have that $D$ covers $W$.
Otherwise, let
$\{d_1, \ldots, d_n\}$ be the set of distribution
nodes of $D$
relative to $W$.  Now we have $V_{\{d_1, \ldots,
d_n\}}(D) \geq
\vert D \vert \geq SN_W(G) \geq \max_{i \in [n]}
V_{\{d_i\}}(W)$,
thus by Theorem~\ref{bigproof} we have that $D$
covers $W$.

Let $q$ be a vertex in $G$ such that
$V_{\{q\}}(W)$ is maximum.
Then the distribution $D$ on $G$ having
$V_{\{q\}}(W) - 1 =
SN_W(G) - 1$ pebbles on $q$ and $0$ pebbles on
all other vertices
fails to cover $W$, as we have $V_{\{q\}}(D) =
V_{\{q\}}(W) - 1 <
V_{\{q\}}(W)$.
\end{proof}

As discussed before, this reduces the problem of
finding
$\Phi_W(G)$ to a matter of computing $$SN_W(G) =
\max_{v \in
V(G)}\sum_{u \in V(G)}W(u)\cdot2^{d(u,v)}.$$ The
corollary that
follows further simplifies the computation of
$\Phi_W(G)$ in
certain cases.  Also, this result provides
encouragement for those
working on a proof for Graham's Conjecture in
standard pebbling.

\begin{definition}
If $W_1$ is a weight function on a graph $G$ and
$W_2$ is a weight
function on a graph $H$, then we define the
weight function $W_1
\times W_2$ on the graph $G \times H$ by
$$[W_1 \times W_2](g,h) =
W_1(g)W_2(h) \quad \forall \, (g,h) \in G \times
H.$$
\end{definition}

\begin{corollary}
Let $G$ and $H$ be connected graphs.  If $W_1$
and $W_2$ are
positive weight functions on $G$ and $H$
respectively, then
$\Phi_{W_1 \times W_2}(G \times H) =
\Phi_{W_1}(G) \Phi_{W_2}(H)$.
\end{corollary}
\begin{proof}
For all $(g,h) \in G \times H$ we have
\begin{eqnarray*}
SN_{W_1 \times W_2}(g,h) &=& \sum_{\substack{g^*
\in
\,G\\h^* \in\, H}}[W_1 \times
W_2](g^*,h^*)2^{d((g,h),(g^*,h^*))}\\
&=& \sum_{h^* \in \, H} \sum_{g^* \in \, G}
W_1(g^*) W_2(h^*)
2^{d(g,g^*)} 2^{d(h,h^*)}\\
&=&\sum_{h^* \in \, H}W_2(h^*)2^{d(h,h^*)}
\sum_{g^* \in\, G}
W_1(g^*)2^{d(g,g^*)}\\
&=&SN_{W_1}(g) \sum_{h^* \in \,
G}W_2(h^*)2^{d(h,h^*)}\\
&=&SN_{W_1}(g) SN_{W_2}(h).
\end{eqnarray*}
The result now follows easily from the stacking
theorem.
\end{proof}

\section{Conjectures}
Our results rely heavily upon the standard value
function, and the
relationship between its values both on $W$ and
$D$ with respect
to particular subsets of $V(G)$.  By
Observation~\ref{value} we
have that $D$ is a cover of $W$ only if
$$ V_S(D) \geq V_S(W) \quad \forall S \subset
V(G).$$
This necessary condition can be proven using only
properties 2 and
3 of $V_S$, as given in Observation~\ref{value}. 
Any function on
the set of distributions having both of these
properties we call a
\emph{general value function}, or simply a value
function.  It can
easily be proven that if $A$ is a value function,
then $D$ covers
$W$ only if
$$ A(D) \geq A(W). $$
We conjecture that the converse is true: that if
$A(D) \geq A(W)$
for all value functions $A$, then $D$ covers $W$.
We leave as an
open question whether a stronger condition is
true: whether if
$V_S(D) \geq V_S(W) \; \forall S \subset V(G)$,
then $D$ covers
$W$.

If $G$ is a connected graph, then $\pi(G)$ can be
thought of as
the maximum of a finite set of values
$\Phi_{W_1}(G), \ldots,
\Phi_{W_n}(G)$. Our results do not apply to the
computation of
these values, as the weight functions involved
are not strictly
positive. However, it is clear that a good
extension of our
results to all weight functions would provide a
determination of
the pebbling numbers of graphs.

\section{Acknowledgements}
This work was done under the supervision of Anant
Godbole at the
East Tennessee State University REU, with
financial support from
the National Science Foundation (Grant
DMS-0139286). The authors
wish to thank Anne Shiu and Nadia Heninger for
their careful
readings of early versions of this paper, and
Anant Godbole for
his guidance and support.

\end{document}